\documentclass[12pt]{article}
\pdfoutput=1

\usepackage{amsmath,amsfonts,amssymb,amsthm,bbm,color,graphicx,url}
\usepackage{enumitem}
\usepackage{natbib}
\usepackage[a4paper]{geometry}
\usepackage{hyperref}
\usepackage{setspace}
\usepackage{booktabs}

\theoremstyle{definition}
\newtheorem{definition}{Definition}[section]
\newtheorem{example}[definition]{Example}
\newtheorem{condition}[definition]{Condition}

\theoremstyle{plain}
\newtheorem{theorem}[definition]{Theorem}
\newtheorem{prop}[definition]{Proposition}

\newcommand{\cF}{\mathcal{F}}

\newcommand{\cI}{\mathcal{I}}
\newcommand{\sA}{\mathsf{A}}
\newcommand{\sO}{\mathsf{O}}
\newcommand{\real}{\mathbb{R}}
\newcommand{\eps}{\varepsilon}
\newcommand{\nat}{\mathbb{N}}
\newcommand{\dd}{\, \mathrm{d}}

\newcommand{\E}{\mathbb{E}}
\renewcommand{\P}{\mathbb{P}}

\newcommand{\GCI}{\mathrm{GCI}}
\newcommand{\SI}{\mathrm{SI}}
\newcommand{\ETI}{\mathrm{ETI}}
\newcommand{\MI}{\mathrm{MI}}
\newcommand{\IS}{\mathrm{IS}}

\newcommand{\one}[1]{\ensuremath{\mathbbm{1}}(#1)}

\DeclareMathOperator{\len}{len}
\DeclareMathOperator{\esssup}{ess\,sup}
\DeclareMathOperator{\essinf}{ess\,inf}

\allowdisplaybreaks

\begin{document}

\begin{center}
\bf \Large 
Scoring Interval Forecasts: Equal-Tailed, Shortest, and Modal Interval

\end{center}

\medskip
\centerline{\bf Jonas R. Brehmer}
\centerline{University of Mannheim, Mannheim, Germany}
\centerline{Heidelberg Institute for Theoretical Studies, Heidelberg, Germany} 

\bigskip

\centerline{\bf Tilmann Gneiting}
\centerline{Heidelberg Institute for Theoretical Studies, Heidelberg, Germany} 
\centerline{Karlsruhe Institute of Technology, Karlsruhe, Germany} 

\bigskip

\centerline{October 23, 2020}

\medskip

\begin{abstract} 

We consider different types of predictive intervals and ask whether
they are elicitable, i.e.\/ are unique minimizers of a loss or scoring
function in expectation.  The equal-tailed interval is elicitable,
with a rich class of suitable loss functions, though subject to
translation invariance, or positive homogeneity and differentiability,
the Winkler interval score becomes a unique choice.  The modal
interval also is elicitable, with a sole consistent scoring function,
up to equivalence.  However, the shortest interval fails to be
elicitable relative to practically relevant classes of distributions.
These results provide guidance in interval forecast evaluation and
support recent choices of performance measures in forecast
competitions.



\end{abstract}

\color{black}

\section{Introduction}  \label{sec:introduction}

In situations where decision making relies on information about
uncertain future quantities, it is desirable to not only have a single
forecast value, i.e.\/ a point forecast, but also information on the
inherent uncertainty of the quantity of interest
\citep{GneitKatz2014}.  A particularly attractive, ubiquitously used
way to achieve this is to require forecasters to report one or
multiple predictive intervals, which are typically designed to contain
the observation with specified nominal probability, as requested
implicitly or explicitly in the Global Energy Forecasting Competition
\citep{Hongetal2016}, the M4 and M5 Competitions \citep{Makrietal2020,
  M5}, and the COVID-19 Forecast Hub \citep{Rayetal2020}.
Consequently, methods for the comparative evaluation of interval
forecasts are in strong demand.  Likewise, researchers and
practitioners need methods for choosing between different models for
the generation of such intervals.

Technically, three types of predictive intervals have been proposed
and used in the literature, two of which are based on the assumption
of a nominal coverage probability $1 - \alpha$, where $\alpha \in
(0,1)$.  The equal-tailed or central interval lies between the
$\frac{\alpha}{2}$- and $(1 - \frac{\alpha}{2})$-quantiles, making it
centered in terms of probability.  The shortest interval has minimal
length, subject to the interval covering the outcome with nominal
probability of at least $1 - \alpha$.  In contrast, the modal interval
maximizes the probability of containing the outcome, subject to a
fixed length.  Early work on the evaluation problem for interval
forecasts can be found in \citet{AitDun1968}, \citet{Winkler1972}, 
\citet{Caselletal1993}, and \citet{Christoff1998}.  Recently, 
\citet{Askanetal2018} have emphasized that tools for the
comparative evaluation of equal-tailed intervals are readily
available, whereas fundamental questions remain open for the shortest
interval.

Consistent scoring (or loss) functions are well-established tools for
quantifying predictive performance and comparing forecasts, see e.g.\/
\citet{DawidMusio2014} and \citet{Gneit2011}
for reviews.  In a nutshell, if we ask forecasters to report a certain
functional of their predictive distributions, then a key requirement
on the loss function is to be (strictly) consistent, in the sense that
the expected loss or score is (uniquely) minimized if the directive
asked for is followed.  The functional is called elicitable if there
is a strictly consistent scoring function.  While consistent scoring
functions have been in routine use for many distributional properties,
such as means or quantiles \citep{Gneit2011}, the existence problem
for any given functional can be a challenge to tackle.  For recent
progress see \citet{Lambetal2008}, \citet{Hein2014}, 
\citet{Steinetal2014}, \citet{FissZieg2016}, and \citet{FronKash2020},
among other works.

The remainder of the paper is structured as follows.
Section~\ref{sec:scoring} provides a short, technical introduction to
the notions of elicitability and consistent scoring functions.  The
core of the paper is in Section~\ref{sec:3intervals}, where we discuss
the elicitability and properties of any consistent scoring functions
for the equal-tailed, shortest, and modal intervals in detail.  We
show that the Winkler interval score arises as a unique choice for the
equal-tailed interval under desirable further conditions, and we
resolve a challenge raised by \citet{Askanetal2018},
who state desiderata for loss functions tailored to the shortest
interval, by showing that in practically relevant settings consistent
scoring functions do not exist.  Although conceptually different, the
modal interval has a close connection to the shortest interval and,
perhaps surprisingly, it has a unique consistent scoring function, up
to equivalence.  Section~\ref{sec:discussion} concludes the paper with
a discussion, where we support the choices of performance measures in
the aforementioned forecast competitions.  Proofs are generally
deferred to the \hyperref[sec:appendix]{Appendix}.

\section{Consistent scoring functions and elicitability}  \label{sec:scoring}

Here we set up notation and provide general technical background and tools.

Let $Y$ be a random variable that takes values in a closed
\textit{observation domain} $\sO \subseteq \real$, and let
$\mathcal{O}$ be the Borel $\sigma$-algebra on $\sO$.  Let $\cF$ be a
class of probability measures on $(\sO, \mathcal{O})$ that represents
the possible distributions for $Y$.  Typically, the observation domain
$\sO$ will either be the real line $\real$, or the set $\nat_0$ of the
nonnegative integers, corresponding to count data, which feature
prominently in applications such as retail and epidemic forecasting.

A statistical property is a functional $T: \cF \to 2^\sA$, where
$2^\sA$ denotes the power set of the \textit{action domain} $\sA
\subseteq \real^n$ that contains all possible reports for $T$.  The
set $T(F) \subseteq \sA$ consists of all correct forecasts for $F \in
\cF$.  Whenever $T(F)$ reduces to a single value $t \in \sA$, we use
the intuitive notation $T(F) = t$ for $T(F) = \{ t \}$.  Moreover, we
let $\E_F$ denote the expectation operator when $Y$ has distribution
$F \in \cF$.  In the special case of an expectation of a derived
binary variable we use the symbol $\P_F$ in customary ways. 
We identify probability distributions with their cumulative
distribution functions (CDFs).

A measurable function $h : \sO \to \real$ is $\cF$-integrable if $\E_F
h(Y)$ is well-defined and finite for all $F \in \cF$.  Finally, a
\textit{scoring function} is a mapping $S : \sA \times \sO \to \real$
such that $S(x, \cdot)$ is $\cF$-integrable for all $x \in \sA$.

\begin{definition}
A scoring function $S$ is \textit{consistent} for a functional $T$
relative to the class $\cF$ if
\begin{equation}   \label{eq:consistency}
\E_F S(t, Y) \le \E_F S(x,Y)  
\end{equation}
for all $F \in \cF$, $t \in T(F)$, and $x \in \sA$.  It is
\textit{strictly consistent} for $T$ if it is consistent for $T$ and
equality in~\eqref{eq:consistency} implies that $x \in T(F)$.  If
there is a scoring function $S$ that is strictly consistent for $T$
relative to $\cF$, then $T$ is called \textit{elicitable}.
\end{definition}

If a forecaster is faced with a penalty $S(x,y)$ for a forecast or
report $x$ and outcome $y$, consistency of the scoring function $S$
for the functional $T$ ensures that any member of the forecaster's 
set of true beliefs $T(F)$ minimizes the expected penalty.  Since the
ordering in~\eqref{eq:consistency} is not affected by scaling $S$ with
a positive constant or adding a report-independent function, we say
that the scoring function $S'$ is \textit{equivalent} to $S$ if
\[
S' (x,y) = c S(x,y) + h(y)
\]
for some $c > 0$ and an $\cF$-integrable function $h : \sO \to \real$.

A basic example of an elicitable functional is the mean or expectation
functional $T (F) := \E_F Y$.  If defined on the class of the
distributions with finite second moment, squared error, $S(x,y) =
(x-y)^2$, is a strictly consistent scoring function for $T$.  As the
mean functional is single-valued, it can be treated in the
point-valued setting, which assumes that functionals map directly into
the action domain $\sA$ \citep{FissZieg2016}.  The following examples
illustrate why interval forecasts call for the full set-valued
framework of \citet{Gneit2011}, which assumes that functionals map
into the power set $2^\sA$.

\begin{example}[quantiles and equal-tailed interval]   \label{ex:quantiles}
For $\alpha \in (0,1)$ an
$\alpha$-quantile of $F$ is a point $x \in \real$ that satisfies $
F(x-) \le \alpha \le F(x)$, where $F(x-) := \lim_{y \uparrow x} F(y)$
denotes the left-hand limit of $F$ at $x$.  The $\alpha$-quantile
functional $T_\alpha (F) := \{ x : F(x-) \le \alpha \le F(x) \}$ is
set-valued, and it is elicitable relative to any class $\cF$.  The
strictly consistent scoring functions are equivalent to
\begin{align}   \label{eq:scsf_quantiles}
S_\alpha (x,y) = \left( \one{y \le x} - \alpha \right) \left( g(x) - g(y) \right),
\end{align}
where $g$ is $\cF$-integrable and strictly increasing, see
\citet{Gneit2011, Gneit2011b} and references therein.  The
equal-tailed interval for $F$ at level $1 - \alpha$ is defined via the
quantiles at levels $\frac{\alpha}{2}$ and $1 - \frac{\alpha}{2}$,
respectively.  Hence, unless both quantiles reduce to single points,
there are multiple equal-tailed intervals at level $1 - \alpha$,
making the equal-tailed interval a set-valued functional, too.
\end{example}

\begin{example}[shortest and modal interval]
Let $\alpha, c \in (0,1)$, and let $F$ be the uniform distribution on
the interval $[0,1]$.  Then every interval of the form $[x, x + 1 -
  \alpha]$, where $x \in [0,\alpha]$, is a shortest interval at level
$1 - \alpha$.  Moreover, every interval of the form $[x, x + c]$,
where $x \in [0, 1 - c]$, is a modal interval at length $c$.
\end{example}

A key characteristic of elicitable functionals is their behavior under
convex combinations of distributions. The following proposition states
the classical convex level sets (CxLS) result
(\citealp[Theorem~6]{Gneit2011}; \citealp{Wang2020}) together with the
refined CxLS$^*$ property of \citet[Proposition~3.3]{Fissetal2019}.

\begin{prop}[convex level sets]  \label{pr:cxls}
Let\/ $T : \cF \to 2^\sA$ be an elicitable functional.  If\/ $F_0, F_1
\in \cF$ and\/ $\lambda \in (0,1)$ are such that\/ $F_\lambda =
\lambda F_1 + (1-\lambda) F_0 \in \cF$, then
\begin{enumerate}[label=(\roman*)]
\item $T(F_0) \cap T(F_1) \subseteq T( F_\lambda)$  (CxLS property);
\item $ T(F_0) \cap T(F_1) \neq \emptyset \,\, \Longrightarrow
  \,\, T(F_0) \cap T(F_1) = T(F_\lambda)$ (CxLS$^{\hspace{0.5mm}*}$
  property).
\end{enumerate}
\end{prop}

If $T$ is a single-valued functional, the properties coincide and are
simply referred to as CxLS.  The most relevant examples of functionals
that do not have convex level sets and thus fail to be elicitable, are
the risk measure Expected Shortfall (ES) and the variance
\citep{Gneit2011}.  If $\sA \subseteq \real$ and certain regularity
conditions hold, convex level sets are also sufficient for
elicitability, as demonstrated by \citet{Steinetal2014}.  However, 
some statistical properties lack
these conditions and fail to be elicitable, even though they have the
CxLS$^*$ property, such as the mode \citep{Hein2014} and tail
functionals \citep{BrehStrok2019}.  In such settings, the following
result can be useful, which is a refined version of Theorem~3.3 of
\citet{BrehStrok2019} that allows for set-valued functionals.

\begin{prop}  \label{pr:Criterion2_sets}
Let\/ $T: \cF \to 2^\sA$ be a functional, and let\/ $F_0, F_1 \in \cF$ be
such that\/ $F_\lambda = \lambda F_1 + (1-\lambda) F_0 \in \cF$ for
all\/ $\lambda \in (0,1)$.  If there are\/ $t_0 \in T(F_0) \backslash
T(F_1)$ and $t_1 \in T(F_1) \backslash T(F_0)$ such that for every\/
$\lambda \in (0,1)$ it holds that either\/ $t_0 \in T(F_\lambda)$
and\/ $t_1 \notin T(F_\lambda)$, or $t_1 \in T(F_\lambda)$ and\/ $t_0
\notin T(F_\lambda)$, then\ $T$ is not elicitable.
\end{prop}

Remarkably, the assertion of Proposition~\ref{pr:Criterion2_sets}
overlaps with part~(ii) of Proposition~\ref{pr:cxls} in the sense that
if $T(F_0) \cap T(F_1) \neq \emptyset$ and the conditions of 
Proposition~\ref{pr:Criterion2_sets} hold, then $T$ cannot have the
CxLS$^*$ property and thus fails to be elicitable. If
$T(F_0) \cap T(F_1) = \emptyset$, Proposition~\ref{pr:Criterion2_sets}
provides a novel result, since Proposition~\ref{pr:cxls}(ii) does not
address this situation.

The criteria for elicitability presented here will be key tools in
what follows.  Like the proofs for the subsequent section, the proof
of Proposition~\ref{pr:Criterion2_sets} is deferred to the
\hyperref[sec:appendix]{Appendix}.

\section{Types of intervals}  \label{sec:3intervals}

We proceed to study equal-tailed, shortest, and modal intervals as
functionals on suitable distribution classes $\cF$.  

Technically, we encode intervals via their lower and upper endpoints
and use the action domain
\[
\sA = \sA_\sO := \{ [a,b] : a, b \in \sO, a \le b \}.  
\]
This choice implies that the predictive intervals we consider are
closed with endpoints in the observation domain $\sO$.  The endpoint
requirement leads to a natural and desirable reduction of the set of
possible intervals for discrete data, such as in the case of count
data, where the endpoints are required to be nonnegative integers.
Closed intervals are compatible with the interpretation of the median
as a `0\% central prediction interval'.  Moreover, in discrete
settings an interval forecast might genuinely collapse to a single
point, so closed intervals allow for a unified treatment of discrete
and continuous distributions.  Lastly, this setting is consistent with
the extant literature, see e.g.\ \citet{Winkler1972},
\citet[Section~7.6]{LambShoh2009}, and \citet{Askanetal2018}.  More
general treatments lead to further complexity without recognizable
benefits.

We denote the length of an interval $I$ as $\len(I)$, and if $\sA'
\subset \sA$ is a set of intervals that all have the same length, we
refer to this common length as $\len (\sA')$.  The left- and
right-hand limits of a function $h: \real \to \real$ at $x$ are
denoted by $h(x-) := \lim_{y \uparrow x} h(y)$ and $h(x+) := \lim_{y
  \downarrow x} h(y)$, respectively.

\subsection{Intervals with coverage guarantees}  \label{sec:GCI} 

A standard principle for interval forecasts is that a correct report
$I$ contains (or covers) the outcome with specified nominal
probability of at least $1 - \alpha$, where $\alpha \in (0,1)$.  A
\textit{guaranteed coverage interval} (GCI) at level $\alpha$ under
the predictive distribution $F$ is any element $[a,b] \in \sA$
satisfying $F(b) - F(a-) \ge 1 - \alpha$, and for all $\eps > 0$
\begin{align}   \label{eq:GCI_restriction}
F(b - \eps) - F(a-) \le 1 - \alpha  
\quad \text{and} \quad    
F(b) - F( (a + \eps)- ) \le 1 - \alpha.
\end{align}
A GCI thus contains just as much probability mass as necessary, but is
not as short as possible. For continuous distributions this definition
reduces to the intuitive requirement $F(b)- F(a) = 1 - \alpha$. We
write $\GCI_\alpha(F)$ for the set of guaranteed coverage intervals at
level $\alpha$ of $F$.  An early theoretical treatment is in
Proposition~7.6 of \citet{LambShoh2009}, according
to which the $\GCI_\alpha$ functional fails to be elicitable relative
to the class of all distributions on the finite domain $\sO = \{ 1,
\ldots, n \}$.  \citet[Section~4.2]{FronKash2019}
apply tools of convex analysis to extend this result to more general
classes of distributions.

It is straightforward to recover these findings by showing that the
$\GCI_\alpha$ functional lacks the CxLS$^*$ property.  Specifically,
let $\alpha \in (0,1)$ and consider continuous distributions $F_0$ and
$F_1$ that satisfy $F_0(b') - F_0(a') = F_1(b') - F_1(a') = 1 - \alpha$ 
for some $a' < b'$, whereas   
\[
F_0 (b) - F_0 (a) > 1 - \alpha \quad \text{and} \quad F_1 (b) - F_1 (a) < 1 - \alpha
\]
for some $a < b$.  Then for some $\lambda \in (0,1)$ we have $[a,b]
\in \GCI_\alpha (F_\lambda)$, even though $[a,b] \not\in
\GCI_\alpha(F_0) \cap \GCI_\alpha(F_1) \not= \emptyset$.  Part (ii) of
Proposition~\ref{pr:cxls} thus implies that the $\GCI_\alpha$
functional fails to be elicitable relative to classes $\cF$ that
contain distributions of the type used here. A similar construction
for discrete distributions is immediate. 

\citet{Fissetal2019} introduce a concept of guaranteed
coverage without the length restriction~\eqref{eq:GCI_restriction},
i.e.\ they consider the class of intervals $[a,b] \in \sA$ which
satisfy $F(b) - F(a-) \ge 1-\alpha$. Like $\GCI_\alpha$, the
corresponding set-valued functional fails to be elicitable
\citep[Corollary~4.7]{Fissetal2019}.

In addition to lacking elicitability, the $\GCI_\alpha$ functional has
the unattractive feature that it fails to be unique for very many
distributions, including, but not limited to, all continuous
distributions.  This motivates the imposition of additional
constraints on the predictive intervals, as discussed now.  \citet{Fissetal2019} discuss still further types of prediction
intervals.

\subsection{Equal-tailed interval (ETI)}  \label{sec:ETI} 

A straightforward way to pick an interval with nominal coverage at
least $1 - \alpha$ under $F$ consists of choosing quantiles at level
$\beta \in (0,\alpha)$ and $\beta + 1 - \alpha$ as the lower and upper
endpoint of the interval, respectively.

The ubiquitous choice is $\beta = \frac{\alpha}{2}$, such that under a
continuous $F$ the outcomes fall above or below the interval with
equal probability of $\frac{\alpha}{2}$.  In general, an
\textit{equal-tailed interval} (ETI) at level $\alpha$ of $F$ is any
member of
\begin{equation}  \label{eq:ETI_def}
\ETI_\alpha (F) := \{ [a,b] \in \sA : a \in T_{\alpha/2}(F), \, b \in T_{1-\alpha/2}(F) \},
\end{equation} 
where $T_\beta (F) := \{ x \in \sO : F(x-) \le \beta \le F(x) \}$
denotes the $\beta$-quantile functional.  The literature also talks of
the `central prediction interval', see, e.g.\ \citet{Fissetal2019}.
In the simplified situation where $F$ is strictly increasing, all
quantiles are unique and thus $\ETI_\alpha(F)$ reduces to a single
interval.

By definition, the $\ETI_\alpha$ functional is equivalent to the
two-dimensional functional $(T_{\alpha/2}, T_{1-\alpha/2} )$, such
that forecasting equal-tailed intervals amounts to forecasting
quantiles.  As a result, the $\ETI_\alpha$ functional is elicitable,
and we can construct consistent scoring functions for it from the
consistent scoring functions~\eqref{eq:scsf_quantiles} for quantiles,
as noted by \citet{GneitRaft2007} and \citet{Askanetal2018}. 
Specifically, if $w_1, w_2$ are
nonnegative weights and $g_1, g_2 : \sO \to \real$ are non-decreasing
$\cF$-integrable functions, then every $S : \sA \times \sO \to \real$
of the form
\begin{align}  \label{eq:ETI_score_general}
S([a,b],y) 
& = w_1 \left( \one{y \le a} - \frac{\alpha}{2} \right) \left( g_1(a) - g_1(y) \right) \\
& \phantom{=} + w_2 \left( \one{y \le b} - \left( 1 - \frac{\alpha}{2} \right) \right) 
  \left( g_2(b) - g_2(y) \right)  \nonumber
\end{align}
is a consistent scoring function for the $\ETI_\alpha$ functional.
Furthermore, $S$ is strictly consistent if $w_1, w_2 \in (0,\infty)$
and $g_1, g_2$ are strictly increasing.  It is no substantial loss of
generality to restrict attention to the class
in~\eqref{eq:ETI_score_general}, since essentially all strictly
consistent scoring functions for $\ETI_\alpha$ are equivalent to this
form.  This is due to the aforementioned fact that $\ETI_\alpha$ can
be interpreted as a vector of two quantiles, and under suitable
regularity conditions, all strictly consistent scoring functions for
vectors of quantiles are equivalent to a sum of scoring functions of
the form~\eqref{eq:scsf_quantiles}, see Proposition~4.2(ii) of \citet{FissZieg2016, FissZieg2019b}.

The choice $w_1 = w_2 = 2/\alpha$ and $g_1(x) = g_2(x) = x$
in~\eqref{eq:ETI_score_general} obtains the classical \textit{interval
  score} (IS) of \citet{Winkler1972}, namely,
\begin{equation}  \label{eq:WIS} 
\IS_\alpha ([a,b], y) := (b-a) 
+ \frac{2}{\alpha} (a - y) \one{y < a} + \frac{2}{\alpha} (y-b) \one{y > b},  
\end{equation} 
which is strictly consistent relative to classes of distributions with
finite first moment.  This is the most commonly used scoring function
for the $\ETI_\alpha$ functional, and scaled or unscaled versions
thereof have been employed implicitly or explicitly in highly visible,
recent forecast competitions \citep{Hongetal2016, Makrietal2020, M5,
  Rayetal2020}.

The Winkler interval score~\eqref{eq:WIS} combines various additional,
desirable properties of scoring functions on $\sO = \real$, such as
\textit{translation invariance}, in the sense that for every $z, y \in
\real$ and $a < b$
\[
S([a-z, b-z], y-z) = S([a,b], y),
\]
and \textit{positive homogeneity} of order~1, in that for every $c >
0$, $y \in \real$, and $a < b$
\[
S( [ca, cb], cy) = c S([a,b], y).
\]
Additionally, the score applies the same penalty terms to values
falling above or below the reported interval, such that it is
\textit{symmetric}, in the sense that
\[
S([a,b], y) = S( [-b, -a], -y)
\]
for $y \in \real$ and $a < b$. 

Our next two results concern scoring functions on $\sO = \real$ that
are of the form~\eqref{eq:ETI_score_general} and share one or more of
these often desirable additional properties.  In particular, the next
theorem demonstrates that either translation invariance or positive
homogeneity and differentiability, combined with symmetry, suffice to
characterize the Winkler interval score~\eqref{eq:WIS}, up to
equivalence.  To facilitate the exposition, assumption (ii) identifies
the action domain $\sA = \{ [a,b] : a \leq b \}$ with the respective
subset $\{ (a,b)' \in \real^2 : a \le b \}$ of the Euclidean plane.

\begin{theorem}  \label{th:IS_characterization}
Let\/ $S$ be of the form~\eqref{eq:ETI_score_general} with
non-constant, non-decreasing functions\/ $g_1$ and\/ $g_2$.  If\/ $S$
is either
\begin{enumerate}[label=(\roman*)]
\item translation invariant, or
\item positively homogeneous and differentiable with respect to\
  $(a,b) \in \sA \subseteq \real^2$, except possibly along the diagonal,
\end{enumerate}
then\/ $g_1$ and\/ $g_2$ are linear. In particular, if\/ $S$ is symmetric
and either (i) or (ii) applies, then\/ $S$ is equivalent to\/
$\IS_\alpha$.
\end{theorem}

The first part of Theorem~\ref{th:IS_characterization}, which states
the linearity of $g_1$ and $g_2$, continues to hold for asymmetric
intervals, defined by choosing endpoints $a \in T_\beta (F)$ and $b \in
T_{\beta + 1 - \alpha}(F)$ for $\beta \in (0,\alpha)$
in~\eqref{eq:ETI_def}.  However, the second statement does not apply,
as non-constant consistent scoring functions for such intervals cannot
be symmetric.

If only symmetry is required in~\eqref{eq:ETI_score_general}, then the
class of possible scoring functions for the equal-tailed interval is
much larger than just the interval score.  To characterize these
functions take $\cI$ to be the class of all non-decreasing functions
$g : \real \to \real$ with the property that $g(x) = \frac{1}{2}
(g(x-) + g(x+))$ for $x \in \real$.  In a trivial deviation from  
\citet{Ehmetal2016} we define the \textit{elementary quantile
  scoring function} as
\[
S_{\alpha, \theta}^\mathrm{Q} (x,y) = ( \one{y \le x} - \alpha ) 
\left( \one{\theta < x} + \frac{1}{2} \one{\theta = x} - \one{\theta < y} 
                        - \frac{1}{2} \one{\theta = y} \right),
\]
which is the special case in \eqref{eq:scsf_quantiles} where $g(x) =
\one{\theta < x} + \frac{1}{2} \one{\theta = x}$.  Given any $\theta
\geq 0$, we now define
\[ 
S_{\alpha,\theta}([a,b], y) = 
S_{\alpha/2, \theta}^\mathrm{Q}(a,y) + S_{1-\alpha/2, - \theta}^\mathrm{Q}(b,y) 
\]
and refer to $S_{\alpha,\theta}$ as the \textit{elementary symmetric
  interval scoring function}.  The following result shows that every
symmetric scoring function of the form~\eqref{eq:ETI_score_general}
arises as a mixture of elementary symmetric interval scoring
functions.  The Winkler interval score \eqref{eq:WIS} emerges in the
special case where the mixing measure $\mu$ is proportional to
Lebesgue measure.

\begin{theorem}  \label{th:IS_symmetry}
Let\/ $S$ be of the form~\eqref{eq:ETI_score_general} with
non-constant, non-decreasing functions\/ $g_1, g_2 \in \cI$.  If\/ $S$ is
symmetric, then it is of the form
\[
S([a,b], y) = \int_{[0,\infty)} S_{\alpha,\theta}([a,b], y) \dd \mu (\theta),
\]
where\/ $\mu$ is a Borel measure on\/ $[0,\infty)$, defined via\/ $\dd
  \mu(\theta) = \dd h(\theta)$ with\/ $h(\theta) = w_1 (g_1(\theta) -
  g_1(-\theta))$ for\/ $\theta \in [0,\infty)$.
\end{theorem}

The usual treatment considers distributions $F \in \cF$ with strictly
increasing CDFs, such that all quantiles are unique.  This ensures
that the interval is truly equal-tailed, with $\ETI_\alpha (F) =
[a,b]$ implying that $\P_F(Y < a) = \P_F(Y > b) = \frac{\alpha}{2}$.
When $F$ admits a Lebesgue density, but some quantiles are not unique,
this property continues to hold.

\begin{table}[t]
\centering
\caption{Properties of the four different intervals in
  $\ETI_\alpha(G)$, where $\alpha = 0.2$.  The expected penalty for an
  interval forecast $[a,b]$ is given by $\E_G \left[ \IS_\alpha
    ([a,b], Y) \one{Y \notin [a,b]} \right]$, so that the expected
  score decomposes into length plus expected penalty.  See text for
  details.  \label{tab:ETI_example}}

\medskip 

\begin{tabular}{ccccc}
\midrule
Interval & Coverage &  Expected $\IS_\alpha$  &  Length &  Expected Penalty \\
\midrule
$[1,2]$  &  0.8  &  3  &  1  &  2  \\
$[0,2]$  &  0.9  &  3  &  2  &  1  \\
$[1,3]$  &  0.9  &  3  &  2  &  1  \\
$[0,3]$  &  1.0  &  3  &  3  &  0  \\
\midrule
\end{tabular}
\end{table}

However, care is needed when interpreting equal-tailed intervals for
discrete distributions.  As a simple example, let $\alpha = 0.2$ and
consider the distribution $G$ on $\nat_0$ that assigns probability
0.1, 0.4, 0.4, and 0.1 to 0, 1, 2, and 3, respectively.  Since
neither the $\frac{\alpha}{2}$- nor the $(1 -
\frac{\alpha}{2})$-quantile are unique, there are four possible
equal-tailed intervals, as listed in Table~\ref{tab:ETI_example}.  The
distribution $G$ illustrates that the coverage of an equal-tailed
interval does not always equal $1 - \alpha$, and may differ among the
valid intervals.  Moreover, $[0,3]$ is not a guaranteed coverage
interval in the sense of Section~\ref{sec:GCI}, as it is unnecessarily
long.  A natural idea is to issue recommendations for such cases,
e.g.\/ `report the shortest available interval' or `report the
interval with the highest coverage'.  However, consistent scoring
functions for the $\ETI_\alpha$ functional cannot be used to ensure
that forecasters follow such further guidelines, since by the
definition of consistency, any valid report attains the same expected
score.

\subsection{Shortest interval (SI)}  \label{sec:SI} 

Instead of defining an interval at the coverage level $1 - \alpha$ via
fixed quantiles, the shortest of these intervals is often sought.
Specifically, a \textit{shortest interval} (SI) at level $\alpha$ of $F$
is any member of the set
\begin{align}  \label{eq:SIargmin}
\SI_\alpha (F) := \arg \underset{[a,b] \in \sA}{\min} \, \{ b-a : F(b) - F(a-) \ge 1 - \alpha \} .
\end{align}
The shortest interval is never longer than an equal-tailed interval,
and in general the two types of intervals differ from each other.  To
see this we follow \citet[Appendix]{Askanetal2018}
and consider a distribution $F$ on $\sO = [0,\infty)$ with strictly
  decreasing Lebesgue density, so that $\SI_\alpha(F) =
  [0,T_{1-\alpha}(F)]$, whereas $\ETI_\alpha (F) = [T_{\alpha/2}(F),
    T_{1-\alpha/2}(F)]$ with a lower endpoint that is strictly
  positive.  However, for distributions with a symmetric, strictly
  unimodal Lebesgue density the two types of intervals are both unique
  and agree with each other.  If a distribution has multiple shortest
  intervals, then neither of them needs to be an equal-tailed
  interval.

As noted in \citet{Askanetal2018}, loss functions
that have been proposed for interval forecasts fail to be strictly
consistent for the $\SI_\alpha$ functional, since they are usually
tailored to the $\ETI_\alpha$ functional.  The question whether the
$\SI_\alpha$ functional is elicitable thus remains unanswered, and
\citet{Askanetal2018} formulate desiderata for
possible scoring functions.  A first result in this direction is
discussed in Section~4.2 of \citet{FronKash2019},
who show that the $\SI_\alpha$ functional fails to be elicitable
relative to classes $\cF$ that contain piecewise uniform
distributions.  In the following we show non-elicitability for more
general classes of distributions, and we also treat discrete
distributions on $\nat_0$.  We start by studying level sets.

\begin{prop}[convex level sets]  \label{prop:SI_CxLS} \hfill 
\begin{enumerate}[label=(\roman*)]
\item The functional\/ $\SI_\alpha$ has the CxLS property.
\item If the class\/ $\cF$ consists of distributions with
  continuous CDFs only, then\ $\SI_\alpha$ has the
  CxLS$^{\hspace{0.5mm}*}$ property.
\end{enumerate}
\end{prop}

The next example shows that the CxLS$^*$ property can be violated for
discrete distributions.

\begin{example}  \label{ex:disc_SI_CxLS}
Let $\alpha \in (0,\frac{1}{3})$, and let $k \geq 1$ be an integer.
Let $\eps \in (0,\frac{\alpha}{3})$ and $\delta \in (0, \eps)$. Let 
$F_0$ and $F_1$ be probability distributions on $\nat_0$ that assign
mass $\eps + \delta$ to $k - 1$ and mass $1 - \alpha - \eps$ to $k$.  
Furthermore,
$F_0$ and $F_1$ assign mass $\eps + \delta$ and $\eps - \delta$,
respectively, to $k + 1$.  This partial specification of $F_0$ and
$F_1$ implies that
\begin{align}    \label{eq:ex:disc_SI_CxLS}
\SI_\alpha(F_0) = \{ [k-1,k] , [k,k+1] \} \quad \text{ and } \quad \SI_\alpha (F_1) = \{ [k-1,k] \},
\end{align}
and for $\lambda \in [0,\frac{1}{2}]$ we have $\SI_\alpha(F_\lambda)
= \SI_\alpha(F_0) \supsetneq \SI_\alpha(F_1)$.  Therefore, 
$\SI_\alpha$ does not have the CxLS$^*$ property relative to any
convex class $\cF$ that includes $F_0$ and $F_1$.
\end{example}

The restrictions on $\alpha$, $\eps$, and $\delta$ in
Example~\ref{ex:disc_SI_CxLS} ensure that the distributions $F_0$ and
$F_1$ are well-defined, unimodal, and 
satisfy~\eqref{eq:ex:disc_SI_CxLS}. To construct such distributions
for general $\alpha \in (0,1)$, we choose $\eps$ and $\delta$ suitably
small and `spread' the probability mass outside of 
$\{ k-1 , k , k+1\}$ such that $k$ is the unique mode 
and~\eqref{eq:ex:disc_SI_CxLS} holds. We thus obtain the following
result.

\begin{theorem}  \label{th:SI_notelic_nat}
Let\/ $k \ge 1$ be an integer, and let\/ $\cF$ be a class of
probability measures on\/ $\nat_0$ that contains all unimodal
distributions with mode\/ $k$.  Then the\/ $\SI_\alpha$ functional is
not elicitable relative to\/ $\cF$.
\end{theorem}

We turn to classes of distributions with Lebesgue densities, so that
the $\SI_\alpha$ functional has the CxLS$^*$ property, and a more
refined analysis proves useful.  First we take up an example in
Section~4.2 of \citet{FronKash2019}.

\begin{example}  \label{ex:uniform}
Given $\alpha \in (0,\frac{3}{5})$, we define distributions $F_0$ and
$F_1$ via the piecewise uniform densities
\[
f_0(x) = (1-\alpha) \mathbbm{1}_{[0,1]}(x) + \frac{\alpha}{3} \mathbbm{1}_{[2,5]}(x)
\;\; \text{and} \;\;
f_1(x) = \frac{1-\alpha}{2} \mathbbm{1}_{[0,2]}(x) + \frac{\alpha}{3} \mathbbm{1}_{[2,5]}(x),    
\]
so that $\SI_\alpha(F_0) = [0,1]$ and $\SI_\alpha(F_1) = [0,2]$,
respectively.  As $\SI_\alpha(F_\lambda) = [0,2]$ for all $\lambda \in
(0,1)$, we conclude from Proposition~\ref{pr:Criterion2_sets} that the
$\SI_\alpha$ functional fails to be elicitable relative to convex
classes of distributions that contain $F_0$ and $F_1$.
\end{example}

As noted, Example~\ref{ex:uniform} applies in situations where the
class $\cF$ includes distributions with piecewise uniform densities.
As this assumption may be restrictive in practice, we proceed to
demonstrate non-elicitability based on substantially more flexible
criteria.

\begin{condition}  \label{cond:technical}
The distribution $F$ admits a Lebesgue density, and there are numbers
$a < b$ and $\eps > 0$ such that $\SI_\alpha(F) = [a,b]$, $F(b) = F(b
+ \eps)$, and if $\beta < \alpha$ then $\len(\SI_\beta(F)) >
\len(\SI_\alpha(F)) + \frac{1}{2} \eps$.
\end{condition}

Loosely speaking, this condition requires that there are `gaps'
on the right- and left-hand side of the shortest interval at
level $\alpha$, while every shortest interval for a level
$\beta < \alpha$ is notably longer than the one at level $\alpha$.

\begin{theorem}  \label{th:SI_notelic_real}
If the class\/ $\cF$ contains the location-scale family of a
distribution satisfying Condition~\ref{cond:technical}, along with its
finite mixtures, then the\/ $\SI_\alpha$ functional is not elicitable
relative to $\cF$.
\end{theorem}

A related result concerning the non-elicitability of $\SI_\alpha$ is
given in Theorem 4.16(i) of \citet{Fissetal2019}.  The
main difference to Theorem~\ref{th:SI_notelic_real} is that 
\citet{Fissetal2019} consider a different class $\cF$ and allow
for scoring functions which take values in the extended real numbers
$\real \cup \{-\infty, \infty\}$.

Although Condition~\ref{cond:technical} might seem technical, suitable
distributions $F$ can be constructed under rather weak assumptions.
For instance, assume $\alpha < \frac{1}{2}$, and let the class $\cF$
contain some compactly supported distribution, along with the
respective location-scale family, and all finite mixtures thereof.
Then constructing an $F$ that satisfies Condition~\ref{cond:technical}
is straightforward.  A more restrictive requirement is the identity
$F(b) = F(b + \eps)$, as it rules out distributions with strictly
positive densities.  The existence of strictly consistent scoring
functions relative to classes of distributions of this type, including
but not limited to the important case of the finite mixture
distributions with Gaussian components, remains an open problem.

We conclude this subsection by considering limit cases of the shortest
interval functional.  For $\alpha \to 0$ the set $\SI_\alpha (F)$
reduces to a single member, namely, the interval $[\essinf(F),
  \esssup(F)]$, where $\essinf$ and $\esssup$ denote the essential
infimum and essential supremum, respectively.  This functional is not
elicitable in our setting \citep{BrehStrok2019}, however, when
allowing for infinite scores, strictly consistent scoring functions
become available \citep[Proposition~4.13]{Fissetal2019}.  For $\alpha
\to 1$ we need to distinguish two cases.  If the elements in $\cF$
admit strictly unimodal densities with respect to Lebesgue measure,
then $\SI_\alpha$ tends to the mode, which fails to be elicitable
\citep{Hein2014}, see also the discussion in Section~\ref{sec:MI}.
For discrete distributions on $\nat_0$ the minimal interval length
zero can be attained so that as $\alpha \to 1$ the members of the set
$\SI_\alpha (F)$ eventually comprise the single point intervals to
which $F$ assigns positive probability.  This limit functional does
not have the CxLS$^*$ property, thus it fails to be elicitable by
Proposition~\ref{pr:cxls}.

\subsection{Modal interval (MI)}  \label{sec:MI} 

In stark contrast to shortest and equal-tailed intervals, we turn to
a type of interval that seeks to maximize coverage, subject to
constraints on length.

Specifically, given any $c > 0$, a \textit{modal interval} (MI) of
length $2c$ of $F$ is any member of the set
\begin{align}  \label{eq:MIargmin}
\MI_c(F) = \arg \underset{[a,b] \in \sA}{\max} \, \{ F(b) - F(a-) : b - a \le 2c \} .
\end{align}
If $F$ has a strictly unimodal Lebesgue density, then the modal
interval shrinks towards the mode as $c \to 0$.  For distributions on
$\nat_0$ the modal interval even agrees with the mode if $c <
\frac{1}{2}$.  This connection and the fact that the length of a modal
interval is fixed, suggest that the $\MI_c$ functional can be
interpreted as a location statistic, whereas the shortest and
equal-tailed intervals contain information on both location and
spread.

In what follows, separate discussions for classes $\cF$ of continuous
and discrete distributions will be warranted.  For distributions on
$\nat_0$, the length of the modal interval will effectively be
$\lfloor 2c \rfloor$, since expanding it further cannot add
probability mass.  In this situation it is convenient to consider $c
\geq 0$, substitute $2c = k$ where $k \in \nat_0$, and encode the
interval via its \textit{lower endpoint} functional $l_k$, so that
$\MI_{k/2}(F) = \{ [x,x+k] : x \in l_k(F) \}$.  Then
\begin{align}  \label{eq:k01_loss}
S(x,y) = - \one{x \le y \le x + k}
\end{align}
is a strictly consistent scoring function for the functional $l_k$ on
the class of all distributions on $\nat_0$.  In particular, the $l_k$
and $\MI_{k/2}$ functionals are elicitable.  In the special case $k =
0$, $l_0$ is the mode functional and~\eqref{eq:k01_loss} becomes
$S(x,y) = - \one{x = y}$, the familiar zero-one or misclassification
loss.  \citet{LambShoh2009} and \citet{Gneit2017} demonstrate that for
distributions with finitely
many outcomes, zero-one loss is essentially the only consistent
scoring function for the mode functional.  We extend this result to
all integers $k \ge 0$, showing that $k$-zero-one-loss
\eqref{eq:k01_loss} is essentially the only strictly consistent
scoring function for the $l_k$ and $\MI_{k/2}$ functionals.

\begin{theorem}   \label{th:MI_scsf_discrete}
Let\/ $k \geq 0$ be an integer, and let\/ $\cF$ be a class of
probability measures on\/ $\nat_0$ that contains all distributions
with finite support.  Then any scoring function that is strictly
consistent for the\/ $l_k$ functional relative to the class\/ $\cF$ is
equivalent to $k$-zero-one-loss\/ \eqref{eq:k01_loss}.
\end{theorem}

For distributions with Lebesgue densities we encode $\MI_c$ via
its \textit{midpoint} functional $m_c$ so that $\MI_c(F) = \{ [x-c,
  x+c] : x \in m_c(F) \}$, where $c > 0$.  Under this convention
\begin{align}  \label{eq:c01_loss}
S(x,y) :=  - \one{ x - c \le y \le x + c}
\end{align}
is a strictly consistent scoring function for $m_c$ on the class of
distributions with Lebesgue densities, whence $m_c$ and $\MI_c$ are
elicitable.  In the limit as $c \to 0$, the scoring
function~\eqref{eq:c01_loss} becomes zero almost everywhere and thus
cannot be strictly consistent for any functional. \citet{Hein2014}
shows that there are no alternative scoring
functions, so the mode fails to be elicitable relative to sufficiently
rich classes of distributions with densities.  Further aspects are
treated in \citet{DearFron2019}.

The following theorem demonstrates, perhaps surprisingly, that
$c$-zero-one-loss \eqref{eq:c01_loss} is essentially the only strictly
consistent scoring function for the $m_c$ and $\MI_c$ functionals.  

\begin{theorem}   \label{th:MI_scsf_continuous}
Let\/ $c > 0$, and let\/ $\cF$ be a class of probability measures on\/
$\real$ that contains all distributions with Lebesgue densities on
bounded support.  Then any scoring function that is strictly
consistent for the\/ $m_c$ functional relative to\/ $\cF$ is almost
everywhere equal to a scoring function which is equivalent to
$c$-zero-one-loss~\eqref{eq:c01_loss}.
\end{theorem}

We complete this section by connecting modal and shortest intervals.
While these are conceptually different types of intervals, a
comparison of~\eqref{eq:SIargmin} and~\eqref{eq:MIargmin} shows that
the $\SI_\alpha$ and $\MI_c$ functionals relate via their defining
optimization problems.  Specifically, the $\SI_\alpha (F)$ functional
is a solution to the constrained optimization problem
\[
\min_{[a,b] \in \sA}  \left( b - a \right) \quad \text{such that} \quad F(b) - F(a-) \ge 1- \alpha,
\]
while the $\MI_c$ functional is a solution to 
\[
\max_{ [a,b] \in \sA} \left( F(b) - F(a-) \right) \quad \text{such that} \quad b - a \le 2c.
\]
Consequently, if either $\len(\SI_\alpha (F)) = 2c$ or $\P_F(Y \in
\MI_c(F)) = 1 - \alpha$, one condition implies the other, and
$\MI_c(F) = \SI_\alpha(F)$ holds.  It remains unclear whether this
connection can be exploited to construct strictly consistent scoring
functions for the $\SI_\alpha$ functional on suitably restrictive,
special classes of distributions.

\section{Discussion}  \label{sec:discussion}

A central task in interval forecasting is the evaluation of competing
forecast methods or models, a problem that is often addressed by using
scoring or loss functions.  For each method or model, and for each
forecast case, the empirical loss is computed.  Losses are then
averaged over forecast cases, and methods with lower mean loss or
score are preferred.  However, for this type of comparative evaluation
to be decision theoretically justifiable, the loss function needs to
be strictly consistent for the predictive interval at hand.

Of the three types of predictive intervals discussed in this paper,
the equal-tailed and modal intervals are elicitable, and we have
discussed the available strictly consistent scoring functions.  For
the popular equal-tailed interval, a rich family of suitable functions
is available, and our findings support the usage of the Winkler
interval score \eqref{eq:WIS}, well in line with implementation
decisions in forecast competitions.  In contrast, the shortest
interval functional fails to be elicitable relative to classes of
distributions of practical relevance.  In this way, we resolve the
questions raised by \citet{Askanetal2018} concerning
the existence of suitable loss functions for the shortest interval in
the negative.  Importantly, there is no obvious way of setting
incentives for forecasters to report their true shortest intervals.
Equal-tailed intervals are preferable due to their elicitability, in
concert with other considerations, such as the intuitive connection to
quantiles and equivariance under strictly monotone transformations
\citep[p.~961]{Askanetal2018}.

The modal interval admits a unique strictly consistent scoring
function relative to comprehensive classes of both discrete and
continuous distributions, up to equivalence.  This appears to be a
rather special situation, as functionals studied in the extant
literature either fail to be elicitable, or admit rich classes of
genuinely distinct consistent scoring functions \citep{FissZieg2016,
  FronKash2019, Gneit2011, Steinetal2014}.  It would be of great
interest to gain an understanding of conditions under which consistent
scoring functions are essentially unique.

As illustrated, interval forecasts are best suited for continuous
distributions, and may exhibit counter-intuitive properties in
discrete settings.  In particular, in the discrete case it may be
unavoidable that the coverage probability of a perfect forecast
exceeds the nominal level $1 - \alpha$.  This raises problems when
assessing interval calibration with the methods of
\citet{Christoff1998}, since asymptotically the null hypothesis of
frequency calibration will then be rejected even under perfectly
correct forecasts.  Modifying the null hypothesis to nominal coverage
greater than or equal to $1 - \alpha$ is not a remedy, since such a
test does not have any power against forecast intervals with too high
coverage.  Consequently, tests for correct forecast specification as
in \citet{Christoff1998} can be problematic when data
fail to be well-approximated by continuous distributions, such as in
the case of retail sales.  Fortunately, comparative evaluation via
consistent scoring functions remains valid and unaffected
\citep{Czadoetal2009, Kolassa2016}.

In many ways, interval forecasts can be seen as an intermediate stage
in the ongoing, transdiciplinary transition from point forecasts to
fully probabilistic or distribution forecasts \citep{Askanetal2018}.
Indeed, probabilistic forecasts in the form of predictive
distributions are the gold standard, as they allow for full-fledged
decision making and well-understood, powerful evaluation methods are
available \citep{Dawid1986, Gneitetal2007, GneitKatz2014}.  Generally,
probabilistic forecasts can be issued in a number of distinct formats,
ranging from the use of parametric distributions, such as in the Bank
of England Inflation Report \citep{Clements2004}, to Monte Carlo
samples from predictive models, as well as simultaneous quantile
forecasts at pre-specified levels, such as in the Global Energy
Forecasting Competition 2014 \citep{Hongetal2016}, the M5 Competition
\citep{M5} and the COVID-19 Forecast Hub \citep{Rayetal2020}.  If the
quantile levels requested are symmetric about the central level of
$\frac{1}{2}$, the collection of quantile forecasts corresponds to a
family of equal-tailed predictive intervals.  Predictive performance
can then be assessed via weighted or unweighted averages of scaled or
unscaled versions of the Winkler interval score~\eqref{eq:WIS}.  The
theoretical results presented here support this widely used practice.

\section*{Appendix: Proofs}  \label{sec:appendix}
\addcontentsline{toc}{section}{Appendix}

\subsection*{Proof of Proposition \ref{pr:Criterion2_sets}}
\addcontentsline{toc}{subsection}{Proof of Proposition \ref{pr:Criterion2_sets}}

Let $t_0, t_1$ be as stated and set $F_\lambda := \lambda F_1 + (1 -
\lambda) F_0$.  Suppose that $S$ is a strictly consistent scoring
function for $T$.  Linearity of expectations in the measure yields
\begin{align*}
\E_{F_\lambda} \left[ S(t_0, Y) - S(t_1, Y) \right] 
& = \lambda \: \E_{F_1} \left[ S(t_0, Y) - S(t_1, Y) \right] \\ 
& \phantom{=} \, + (1 - \lambda) \: \E_{F_0} \left[ S(t_0, Y) - S(t_1, Y) \right],
\end{align*}
where the first difference is positive, while the second difference is
negative.  Consequently, $\E_{F_\lambda} S(t_0, Y) = \E_{F_\lambda}
S(t_1, Y)$ for some $\lambda \in (0,1)$.  Since either $t_0 \in
T(F_\lambda)$ and $t_1 \notin T(F_\lambda)$, or $t_1 \in T(F_\lambda)$
and $t_0 \notin T(F_\lambda)$, we arrive at a contradiction.

\subsection*{Proof of Theorem \ref{th:IS_characterization}}
\addcontentsline{toc}{subsection}{Proof of Theorem \ref{th:IS_characterization}}

Let $S$ be a scoring function of the form
\eqref{eq:ETI_score_general}.  Let $y, z \in \real$, $a < b$ and
choose $b = y$.  Then translation invariance of $S$ gives
\begin{align*}
- w_1 \frac \alpha 2 (g_1(a) - g_1(y)) 
& = S([a,y], y) \\
& = S([a-z, y-z], y-z) \\
& = - w_1 \frac \alpha 2 (g_1(a-z) - g_1(y-z)),
\end{align*}
and rearranging yields $g_1(a) - g_1(y) = g_1(a-z) - g_1(y-z)$ for $a,
y, z \in \real$.  Choose $y = 0$ and define $\tilde{g}(x) := g_1(x) -
g_1(0)$ to obtain $\tilde{g}(a-z) = \tilde{g}(a) + \tilde{g}(-z)$ for
$a, z \in \real$.  Thus $\tilde{g}$ obeys Cauchy's functional
equation, and since $\tilde{g}$ is non-constant and non-decreasing, we
get $g_1(x) = \gamma x + g_1(0)$ for some $\gamma > 0$.  For $g_2$ we
apply the same arguments, to complete the proof of part (i).

Let $y \in \real$, $a < b$ and choose $b = y$.  If $S$ is positively
homogeneous then for all $c > 0$
\begin{align*}
- w_1 c \frac \alpha 2 (g_1(a) - g_1(y))  
& = c S([a,y], y) \\
& = S([ca, cy], cy) = - w_1 \frac{\alpha}{2} (g_1 (ca) - g_1(cy)),
\end{align*}
and thus $c (g_1(a) - g_1(y)) = g_1(ca) - g_1(cy)$ for $a, y \in
\real$ and $c > 0$.  Choose $y = 0$ and define $\tilde{g}(x) := g_1(x)
- g_1(0)$ to obtain $c \tilde{g}(a) = \tilde{g}(ca)$ for $c > 0$ and
$a \in \real$, as in Section C of the Supplementary Material for \citet{NoldeZieg2017}.  Since $\tilde{g}$ is non-constant,
non-decreasing, and differentiable, $g_1(x) = \gamma x + g_1 (0)$ for
some $\gamma > 0$.  Using the same arguments for $g_2$ we complete the
proof of part (ii).

Now suppose $S$ is also symmetric and $g_2(x) = \rho x + g_2(0)$ for
some $\rho > 0$.  Then the same reasoning as in the proof of
Theorem~\ref{th:IS_symmetry} shows that $w_1 \gamma = w_2 \rho$, which
proves the equivalence to $\IS_\alpha$.

\subsection*{Proof of Theorem \ref{th:IS_symmetry}}
\addcontentsline{toc}{subsection}{Proof of Theorem \ref{th:IS_symmetry}}

Let $S$ be a scoring function of the form~\eqref{eq:ETI_score_general}
and let $a,b, y \in \real$ with $a<b$ and $b = y$.  Then the symmetry
of $S$ gives
\begin{align*}
- w_1 \frac \alpha2  (g_1 (a) - g_1 (y) ) &= S ( [a,y], y) \\
&= S( [-y, -a], -y) = w_2 \frac \alpha2 ( g_2(-a) - g_2 (-y) ) ,
\end{align*}
and rearranging yields $w_1 (g_1 (a) - g_1(y)) = w_2 (g_2(-y) -
g_2(-a))$ for $a, y \in \real$.  For $x, y, \theta \in \real$, define
the function
\begin{align*}
f (x,y, \theta) := \one{\theta < x} + \frac{1}{2} \one{\theta = x} 
                 - \one{\theta < y} - \frac{1}{2} \one{\theta = y},
\end{align*}
which satisfies $f(-y,-x,\theta) = f(x,y,-\theta)$ for $x, y, \theta
\in \real$. Recall that $\cI$ is the class of non-decreasing functions
$g : \real \to \real$ such that $g(x) = \frac{1}{2} (g(x-) + g(x+))$
for $x \in \real$. For all $g \in \cI$ and $y < x$
\begin{align*}
\int f(x,y,\theta) \dd \mu_g (\theta) 
= \frac{1}{2} (g(x+) - g(y-)) + \frac{1}{2} (g(x-) - g(y+)) 
= g(x) - g(y),
\end{align*}
where $\mu_g$ is the Borel measure on $\real$ induced by $g$.  If we
define the measures $\mu_1 = w_1 \mu_{g_1}$ and $\mu_2 = w_2
\mu_{g_2}$, then the first part of the proof implies
\begin{align*}
\int f(x,y, \theta ) \dd \mu_2 (\theta) &= w_2 ( g_2 (x) - g_2(y) )
\\ &= w_1 ( g_1 (-y) - g_1 (-x) ) \\ &= \int f(-y, -x, \theta ) \dd
\mu_1 (\theta) = \int f(x , y, - \theta ) \dd \mu_1 (\theta )
\end{align*}
for $y < x$, and the proof is completed by defining $\mu$ via
$\mu((y,x]) = \mu_1((y,x]) + \mu_1([-x,-y))$.

\subsection*{Proof of Proposition~\ref{prop:SI_CxLS}}
\addcontentsline{toc}{subsection}{Proof of Proposition~\ref{prop:SI_CxLS}}

Let $F_0, F_1 \in \cF$, and suppose that $[a,b] \in \SI_\alpha(F_0)
\cap \SI_\alpha(F_1)$.  Set $F_\lambda := \lambda F_1 + (1 - \lambda)
F_0$ and note that for all $\lambda \in (0,1)$ and all $s, t \in
\real$ we have
\begin{align}   \label{eq:FconvexCombi}
F_\lambda (t) - F_\lambda (s-) = \lambda \left( F_1 (t) - F_1 (s-) \right) 
+ (1 - \lambda) \left( F_0 (t) - F_0 (s-) \right).
\end{align}
In particular, $F_\lambda (b) - F_\lambda (a-) \geq 1 - \alpha$ and
$[a,b] \in \SI_\alpha(F_\lambda)$, as otherwise \eqref{eq:FconvexCombi}
yields a contradiction to our initial assumption.  This proves part
(i).

Now let $F_0, F_1 \in \cF$ have continuous CDFs. Since $(s,t) \mapsto
F(t) -F(s)$ is a continuous function for all $F \in \cF$, we must have
$F (b) - F (a) = 1-\alpha$ for every $[a,b] \in \SI_\alpha (F)$.
Suppose $[a',b'] \in \SI_\alpha(F_0) \cap \SI_\alpha(F_1)$, as
otherwise there is nothing to show, and let  $\lambda \in (0,1)$ and
$[a,b] \in \SI_\alpha(F_\lambda)$ be given. By the first part of the
proof
\begin{equation}  \label{eq:lenSI_contCDF}
\len(\SI_\alpha(F_1)) = \len(\SI_\alpha(F_0)) = b' - a' = b - a.
\end{equation}
Furthermore, $F_\lambda (b) - F_\lambda (a) = 1-\alpha$ and we see
from~\eqref{eq:FconvexCombi} that $F_0 (b) - F_0 (a) \ge 1-\alpha$ or
$F_1(b) - F_1(a) \ge 1 - \alpha$ must hold. Suppose the first of these
two inequalities is satisfied. Then equality must hold since the
strict inequality $F_0(b) - F_0(a) > 1-\alpha$ would 
contradict~\eqref{eq:lenSI_contCDF}. This yields $[a,b] \in
\SI_\alpha (F_0)$ and via~\eqref{eq:FconvexCombi} we obtain $F_1 (b)
- F_1(a) = 1-\alpha$. Taken together this gives $[a,b] \in
\SI_\alpha(F_1)  \cap \SI_\alpha(F_0)$, which proves part (ii).

\subsection*{Proof of Theorem \ref{th:SI_notelic_real}}
\addcontentsline{toc}{subsection}{Proof of Theorem \ref{th:SI_notelic_real}}

We proceed by constructing suitable convex combinations as in
Example~\ref{ex:uniform}.  Specifically, let $F_0$ satisfy
Condition~\ref{cond:technical}, and without loss of generality assume
that $\SI_\alpha(F_0) = [0,b]$ for some $b > 0$.  For instance, if
$\alpha \in (0, \frac 12)$ a valid choice is $F_0 = \alpha G_1 +
(1-\alpha) G_0$, where $G_0$ and $G_1$ are absolutely continuous
distributions with support $[0,b]$ and $[2b,3b]$,
respectively. Define $F_1$ via
\[
F_1(x) := F_0 \left( \frac{b}{b + \frac{1}{2} \eps} \, x \right)
\]
and set $F_\lambda := \lambda F_1 + (1 - \lambda) F_0$.  We proceed to
show that $[0, b + \frac{1}{2} \eps] \in \SI_\alpha(F_\lambda)$ for
all $\lambda \in (0,1]$, which allows us to apply
  Proposition~\ref{pr:Criterion2_sets} and conclude non-elicitability.

Clearly, $\SI_\alpha(F_1) = [0, b + \frac{1}{2} \eps]$, and since
$F_0(b) = F_0(b + \eps)$ it holds that $F_\lambda(b + \frac{1}{2}
\eps) - F_\lambda(0) = 1 - \alpha$ for $\lambda \in (0,1)$.  For a
contradiction, suppose there are $\lambda \in (0,1)$ and $a_\lambda
\leq b_\lambda$ with $F_\lambda(b_\lambda) - F_\lambda(a_\lambda) \geq
1 - \alpha$ and $b_\lambda - a_\lambda < b + \frac{1}{2} \eps$.  Since
$\SI_\alpha(F_1) = [0, b + \frac{1}{2} \eps]$ it cannot be true that
$F_1(b_\lambda) - F_1(a_\lambda) \ge 1 - \alpha$ and so
$F_0(b_\lambda) - F_0(a_\lambda) > 1 - \alpha$ must hold, for a
contradiction to the final part of Condition~\ref{cond:technical}.
Consequently, $\SI_\alpha (F_\lambda) = [0, b + \frac{1}{2} \eps]$ for
all $\lambda \in (0,1]$, and the proof is complete.

\subsection*{Proof of Theorem \ref{th:MI_scsf_discrete}}
\addcontentsline{toc}{subsection}{Proof of Theorem \ref{th:MI_scsf_discrete}}

Let $k \geq 0$ be an integer, and suppose that $S$ is a strictly
consistent scoring function for the functional $l_k$ relative to
$\cF$.  To facilitate the presentation, we introduce the alternative
notation $S(M,y)$ for $S(x_M,y)$, where $x_M \in \nat_0$ denotes the
lower endpoint of an interval $M \in \sA $, with
$\sA = \{ [x, x+k] : x \in \nat_0 \}$.  We proceed in three steps.

\paragraph{Step 1}  We show that $S$ is of the form 
\begin{equation}  \label{eq:formofS}
S(x,y) = g(x,y) \one{x \le y \le x + k} + h(y)
\end{equation}
for functions $g : \nat_0 \times \nat_0 \to \real$ and $h : \nat_0 \to
\real$.

To this end, let $M_0, M_1 \in \sA$ such that $M_0 \cap M_1 =
\emptyset$.  For a contradiction, suppose that the mapping $\varphi :
\nat_0 \to \real$ defined via $\varphi (y) = S(M_0, y) - S(M_1, y)$ is
non-zero on $U := (M_0 \cup M_1)^c \cap \nat_0$.  We first treat the
case where
$\varphi (y) = c$ for all $y \in U$ and some $c \in \real \backslash
\{0 \}$.  If $c > 0$ let $F_0$ be the uniform distribution on $M_0$
and for all $n \in \nat$ let $F_n$ be the uniform distribution on
some set $U_n \subset U$ with $\vert U_n \vert = 2nk$. If we define
$G_n := \frac 1n F_0 + (1- \frac 1n) F_n$, then $\MI_{k/2} (G_n) =
\MI_{k/2} (F_0) = M_0$ for all $n \in \nat$.  Since $\int \varphi (y)
\dd G_n (y) \to c > 0$ for $n \to \infty$, we obtain a contradiction
to the strict consistency of $S$.  A similar argument applies if $c <
0$. Consequently, $\varphi$ cannot be constant on $U$, i.e.\ there are
$i_0, i_1 \in U$ such that $\varphi (i_0) \neq \varphi (i_1)$.

Now set $I := \{i_0, i_1\}$. As the class $\cF$
contains all distributions with finite support, we
can find probability measures $F_0, F_0', F_1 \in \cF$ that
satisfy the following three conditions:
\begin{enumerate}[label=(\roman*)]
\item There exists a $\lambda^* \in (0,1)$ such that for
  $F_\lambda := \lambda F_1 + (1-\lambda) F_0$ and $F_\lambda' :=
  \lambda F_1 + (1-\lambda) F_0'$
  \[
  \MI_{k/2}(F_\lambda) = \MI_{k/2}(F_\lambda') = 
  \left\{ \begin{array}{ll} M_0, & \lambda < \lambda^*, \\ 
                            M_1, & \lambda > \lambda^*. \end{array} \right.
  \]
\item $F_0$ and $F_0'$ coincide outside of $I$.
\item $\int_I \varphi (y) \dd F_0 (y) \neq \int_I \varphi (y) \dd F_0' (y)$.
\end{enumerate}

To see this, define $F_0$ and $F_0'$ via the probabilities $F_0 (\{j\})
= F_0' (\{j\}) = 1 / (k+2)$ for $j\in M_0 $ and
\begin{align}  \label{eq:F0specification}
F_0 (\{i_0\}) = F_0'(\{i_1\}) 
= \frac{1}{2(k+2)} + \eps,  \quad \text{and} \quad F_0 (\{i_1\}) 
= F_0'(\{i_0\}) = \frac{1}{2(k+2)} - \eps 
\end{align}
for some $\eps \in (0,1/( 2 (k+2)) )$. Condition~(ii) is immediate and
(iii) follows from the fact that $\varphi (i_0) \neq \varphi (i_1)$.
Moreover, letting $F_1$ be the uniform distribution on $M_1$ ensures
(i).

Consider the integrated score difference
\[
\Delta (F, G, \lambda) := \int \left( S(M_0, y) - S(M_1, y) \right) 
\dd (\lambda G + (1-\lambda)F)(y),
\]
which is linear in $\lambda \in [0,1]$.  The strict consistency of $S$
in concert with (i) yields $\Delta (F_0, F_1, 0) < 0$, $\Delta( F_0',
F_1, 0) < 0$, and $\Delta (F_0, F_1, 1) = \Delta (F_0', F_1, 1) > 0$.
Since $\Delta (F_0, F_1, 0)$ $\neq \Delta(F_0', F_1, 0)$ by (ii) and
(iii), the linear mappings $\lambda \mapsto \Delta (F_0, F_1,
\lambda)$ and $\lambda \mapsto \Delta (F_0', F_1, \lambda)$ must have
distinct roots.  This implies that one of the two mappings does not
vanish at $\lambda^*$, in contradiction to the consistency of $S$.
Consequently, $\varphi = 0$ on $U$ such that we can conclude $S(M_0,
y) = S(M_1, y)$ for all $y \in (M_0 \cup M_1)^c$.  By varying the
disjoint intervals $M_0, M_1 \in \sA$, we obtain that for all $y \in
\nat_0$ the values $S(M,y)$ are the same for all $M \in \sA$ with $ y
\notin M$. This yields that there exists a function $h: \nat_0 \to
\real$ such that $S$ is of the form~\eqref{eq:formofS}.

\paragraph{Step 2}  Now we prove that $y \mapsto g(x,y)$ is constant on
$[x, x+k]$.  As before, we use the notation $g(M,y)$ for $g(x_M,y)$,
where $x_M \in \nat_0$ is the lower endpoint of $M \in \sA$. For $k =
0$ there is nothing to show, so let $k> 0$. For a contradiction,
suppose there is an $M_0 \in \sA$ such that $y \mapsto g(M_0, y)$ is
not constant on $M_0$, i.e.\ there are $i_2, i_3 \in M_0$ such that
$g(M_0, i_2) \neq g(M_0, i_3)$.  This ensures that we can choose an
interval $M_1 \in \sA$, with $M_1 \cap M_0 = \emptyset$, and
distributions $F_0, F_0', F_1 \in \cF$ that satisfy conditions~(i),
(ii), and~(iii) in Step~1, for $I = \{ i_2, i_3 \}$.  For example, we
can choose $F_0$ and $F_0'$ by using the uniform distribution on $M_0$
and modifying it at $i_2$ and $i_3$ as in~\eqref{eq:F0specification},
while ensuring $M_1$ is separated from $M_0$ by a sufficiently large
gap.  As in Step~1 we obtain $\Delta(F_0, F_1, 0) \neq \Delta(F_0',
F_1, 0)$ such that the mappings $\lambda \mapsto \Delta (F_0, F_1,
\lambda)$ and $\lambda \mapsto \Delta (F_0', F_1, \lambda)$ have
distinct roots.  This is a contradiction to the consistency of $S$ and
proves that $y \mapsto g(M_0, y)$ is constant on $M_0$.  We can thus
replace $g(x,y)$ in~\eqref{eq:formofS} by $\tilde{g}(x)$ for some
function $\tilde{g} : \nat_0 \to \real$.

\paragraph{Step 3}  It remains to be shown that $\tilde{g}$ reduces 
to a negative constant.  To this end, consider $M_0 \in \sA$ and $M_1
\in \sA$ and assume that $\tilde{g}(M_0) < \tilde{g}(M_1)$.  Due to
the specific form of \eqref{eq:formofS} we have
\begin{align*}
\E_F \left[ S(M_0, Y) - S(M_1,Y) \right]  
= \tilde{g}(M_0) \P_F(Y \in M_0) - \tilde{g}(M_1) \P_F(Y \in M_1)
\end{align*}
for all $F \in \cF$.  However, due to the strict consistency of $S$
this expression must be negative if $M_0 \in \MI_{k/2}(F)$ and positive
if $M_1 \in \MI_{k/2}(F)$, for the desired contradiction.  Therefore
$\tilde{g}$ reduces to a constant, and using once more the consistency
of $S$, we see that this constant is negative.  The proof is complete.

\subsection*{Proof of Theorem \ref{th:MI_scsf_continuous}}
\addcontentsline{toc}{subsection}{Proof of Theorem \ref{th:MI_scsf_continuous}}

We sketch this proof only, as it proceeds in the very same three steps
as the proof of Theorem \ref{th:MI_scsf_discrete}.  Specifically, let
$c > 0$, and let $S$ be a strictly consistent scoring function for the
functional $m_c$ relative to $\cF$.  In Step 1, we show that $S$ is 
almost everywhere of the form
\[
S(x,y) = g(x,y) \one{x-c \le y \le x+c} + h(y)
\]
for $\cF$-integrable functions $g : \real \times \real \to \real$ and
$h : \real \to \real$.  In Step 2 we prove that $g$ reduces to a
function $\tilde{g}$ in the variable $x$ only, and in Step 3 we
demonstrate that $\tilde{g}$ reduces to a negative constant.  The
technical details are analogous to those in the above proof of 
Theorem~\ref{th:MI_scsf_discrete}, with the only difference that
the set $I$ is now an interval and the statements hold Lebesgue
almost everywhere.

\section*{Acknowledgments}
\addcontentsline{toc}{section}{Acknowledgments}

The authors are grateful for support by the Klaus Tschira Foundation.
Jonas Brehmer gratefully acknowledges support by the German Research
Foundation (DFG) through Research Training Group RTG 1953.  We thank
two anonymous referees, Francis Diebold, and Tobias Fissler for
thoughtful comments and suggestions.

\end{document}